\input AHTOHFIE.STY

\UDC{%
515.142.331
+515.142.321
+512.543.16
%
%
}


\title{%
Forester's lattices
and
small non-Leighton complexes
}
\author{%
Natalia S. Dergacheva
\quad
Anton A. Klyachko
}
\address{%
Faculty of Mechanics and Mathematics of Moscow State University
\\
Moscow 119991, Leninskie gory, MSU.
\\
Moscow Center for Fundamental and Applied Mathematics
\\
nataliya.dergacheva@gmail.com
\quad
klyachko@mech.math.msu.su
}

\grants{\RSF 22-11-00075}

\abstract{%
We construct two CW-complexes $K$ and $L$ admitting a common, but not
finite common, covering, where $K$ is homeomorphic to a complex with a
single 2-cell.
}

\s 1.
Introduction

\proclaim Leighton's theorem {\rm[Lei82]}.
If two finite graphs have a common covering, then they have a common
finite covering.

\noindent
Alternative proofs and generalisations of this result can be
found, e.g., in
[Neu10],
[BaK90],
[SGW19],
[Woo21],
[BrS21],
and references therein.

\noindent
For two-dimensional CW-complexes
(and cell coverings)
the similar assertion is false:
\-
the first example of
such a \emph{non-Leighton
pair} of complexes
(i.e., a pair of finite CW-complexes
having a common, but not finite common, covering)
contained six two-cells
in each complex ([Wis96] and [Wis07]);
\-
later, this number was reduced to four
[JaW09];
\-
and then to two [DK23], [DK24].

\enditem
It is unknown
whether there exists a non-Leighton pair
with one two-cell (in both complexes
or at least in one of them).
We give an ``almost answer" to this open question from [DK23]:
\disp{\sl
there exists a non-Leighton pair $(K,L)$,
in which the complex $K$
is
\spacing{homeomorphic} to a complex with one two-cell.
}%
In addition, our example has the following
extremal properties:
\-
the single-two-cell complex $\=K$ homeomorphic to $K$
is the standard (one-vertex) complex of the
presentation $\BS(2,4)\:=\pres<c,d|d^{-1}c^2d=c^4>$,
so $\=K$ (as it is easy to see)
is
\emph{Leighton},
i.e., it cannot be a member of any non-Leighton pair;
the non-Leighton complex $K$
is obtained from $\=K$ by adding
an edge dividing the 2-cell into two parts;
this is the first (and unexpected) example, where a
cell subdivision
causes such a dramatic effect;
\-
the complex $K$ contains two 2-cells
(this is the minimum known example, see above;
however, the complex $L$ is fairly large,
while, in [DK23] (and [DK24]),
a non-Leighton pair, in which both complexes
have two 2-cells, was constructed);
\-
finally,
the fundamental group of $K$ is one-relator,
and the fundamental group of $L$ is commensurable with
a one-relator group, i.e.,
passing to a suitable
finite covering of $L$, we obtain
the following nontrivial fact:
\dispno{\sl\narrower\narrower
there exists a non-Leighton pair of complexes $(K,\~L)$,
in which $\pi_1(K)$ is a one-relator group,
and $\pi_1(\~L)$ is a finite-index subgroup in a
one-relator group.
}(1)%
Compare this with the Bridson--Shepherd theorem [BrS22]
(the version from [DK24]):
\disp{\sl\narrower\narrower
in a non-Leighton pair, the fundamental group of neither complex
can be free \(or even virtually free\).
}%
The following natural question remains open.

\Question.
Does there exist a non-Leighton pair, in which the fundamental groups
of both complexes are one-relator?

Fact (1) itself does not answer this question: e.g., the free
products $(\Z/3\Z)*(\Z/3\Z)$ and $(\Z/5\Z)*(\Z/5\Z)$ cannot be
realised as fundamental groups of finite complexes with a common
covering [MSW03], although their finite-index subgroups can (these free
products are even commensurable).

\smallskip

Actually our example of the
non-Leighton pair $(K,L)$ follows easily from the
arguments of
Forester's paper [Fo24].
Our humble contribution comprises
\-
noticing this fruit of Max Forester's work [Fo24]
(he did not seem interested in minimal non-Leighton
pairs)
\-
and writing down a short and
almost self-contained proof, where
``almost self-containedness"
means that we use only
standard facts and
the following special case
of the classification of the Baumslag--Solitar groups up to
commensurability [CKZ21] (see also [Fo24]):
\disp{\sl
the Baumslag--Solitar groups $\BS(2,4)$ and $\BS(4,16)$
are \emph{incommensurable}
\(i.e., no finite-index subgroup of the first group is
isomorphic to a finite-index subgroup of the latter one\).}

{\noindent \bf Our notation}
is mainly standard. Note only that, if
$k\in \Z$ and $x$ and $y$ are elements of a group, then $x^y$,
$x^{ky}$, and $x^{-y}$ denote $y^{-1}xy$, $y^{-1}x^ky$, and
$y^{-1}x^{-1}y$, respectively.
A commutator $[x,y]$ denotes $x^{-1}y^{-1}xy$.
The \emph{Baumslag--Solitar groups} are
$\BS(n,m)\:=\pres<c,d|c^{nd}=c^m>$.

We are grateful to an anonymous referee for
valuable remarks allowing us to improve the readability of the paper
and correct some errors.
The authors thank the
Theoretical Physics and Mathematics Advancement Foundation ``BASIS".

\s 2.
Main theorem

Consider the following 2-complexes $K$ and $L$.

The complex $K$ is the standard (one-vertex) complex
for the presentation
$\pres<c,d,z|c^{-1}dc^2=z=cdc^{-2}>$
(in [Fo24], this complex was denoted $Z_{2,4}$).
Clearly, this two-2-cell complex is homeomorphic
to a one-2-cell one:
the standard complex of the standard presentation
$\BS(2,4)\:=\pres<c,d|c^{2d}=c^4>$.

\noindent
The complex $L$ consists of
\-
two vertices: \emph{black} and \emph{white};
\-
six edges:
\itemitem{--}
the edges $c_\bullet $ and $c_\circ $ are loops at
the black and white vertices,
respectively;
\itemitem{--}
the edges $y$ and $z$ start at the black vertex
and end at the white one;
\itemitem{--}
the edges $t$ and $t_1$ start at the white vertex and end
at the black one;
\-
four two-dimensional cells $A,B,C,D$ with boundaries
$y^{-1}c_\bullet yc_\circ ^{-2}$,\quad
$z^{-1}c_\bullet zc_\circ ^{-2}$,\quad
$t_1^{-1}c_\circ tc_\bullet ^{-2}$,\quad
$t^{-1}c_\circ t_1c_\bullet ^{-2}$
(Fig. 1).

\goodbreak
\bigskip
\centerline{\input 1.PIC}
\nobreak%
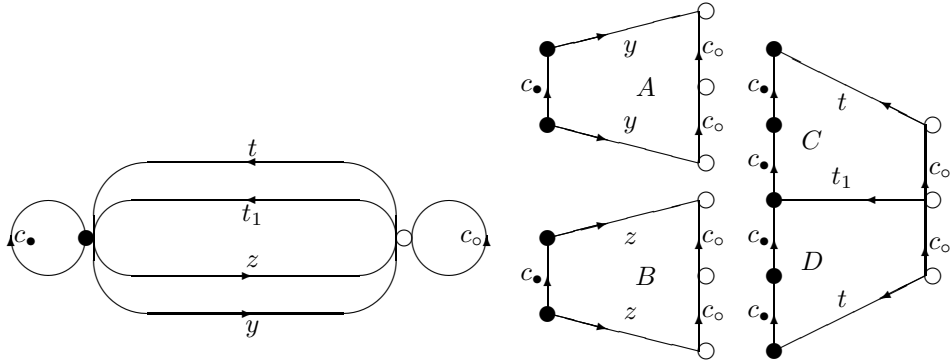
\centerline{Fig. 1}%
\goodbreak
\bigskip

\enditem
(In [Fo24], this complex was nameless,
it is the quotient complex of~$X_{2,4}$
by the action of the group $G_2$.)

\Theorem.
The complexes $K$ and $L$ form a \emph{non-Leighton pair}\/:
\-
their universal coverings are isomorphic
\-
and they have no finite common covering.
\enditem
Moreover, the complex $K$ containing two 2-cells
is homeomorphic to a complex with one 2-cell
\(and the complex $L$ containing four 2-cells is
homeomorphic to a complex with three 2-cells\).

The next two sections constitute the proof of this theorem.

\s 3.
Why the complexes $K$ and $L$ have no common finite covering?

Because their fundamental groups are incommensurable.
The fundamental group of $K$ is $\BS(2,4)$, surely;
while the fundamental group of $L$ has the following presentation:
$$
\pi_1(L)=
\pres<c_\bullet ,c_\circ ,y,z,t,t_1|
y^{-1}c_\bullet yc_\circ ^{-2}=1,\;
z^{-1}c_\bullet zc_\circ ^{-2}=1,\;
t^{-1}c_\circ t_1c_\bullet ^{-2}=1,\;
t_1^{-1}c_\circ tc_\bullet ^{-2}=1,\;
z=1
>.
$$
Obvious simplifications
(Tietze's transformations) give the following presentation:
$$
\pi_1(L)=
\pres<c_\circ ,y,t|
[y,c_\circ ^2]=1,\;
c_\circ ^{2t}=c_\circ ^8
>.
$$
The kernel $H$ of the homomorphism from $\pi_1(L)$ onto $\Z_2\:=\Z/2\Z$
mapping $c_\circ $ to 1, and $y$ and $t$ to 0
has index two in~$\pi_1(L)$.
If we take the set $T=\{1,c_\circ^{-1}\}$ as a Schreier transversal,
then we obtain the following generating set
for $H$:
$$
\bigl\{\tau x(\={\tau x})^{-1}\;|\;x\in\{c_\circ,y,t\}\bigr\}
\setminus\1
=
\left\{q=c_\circ^2,\;y,\;t,\;\^y=y^{c_\circ},\;\^t=t^{c_\circ}\right\}
$$
(where $\=v\in T$ is the only element of $T\cap Hv$).
The defining relators for
$H$ are
words of the form $\tau r\tau^{-1}$
(expressed via the generators of $H$), where $r$ is a defining relator
of the initial group and
$\tau\in T$ [Sch27] (See also [LS77]). Thus, we obtain
$$
H=
\pres<q,y,t,\^y,\^t|[y,q]=1,\;[\^y,q]=1,\;q^t=q^4,\;q^{\^t}=q^4>.
$$
Tietze transformations allows us to simplify this presentation:
$$
\eqalign{
H&\iso
\pres<q,y,t,\^y,\^t, y',y''|[y,q]\=^3 1\=^4 [\^y,q],\;q^t=q^{\^t}=q^4,\;
y'=yt,\;y''=\^yt,\;q^{y'}\=^1 q^4\=^2 q^{y''}
>
\iso
\cr
\iso&
\pres<q,y,t,\^y,\^t,y',y''|q^t=q^{\^t}=q^4,\;
y'\=^5 yt,\;y''\=^6 \^yt,\;q^{y'}=q^4=q^{y''}
>
\iso
\pres<q,t,\^t,y',y''|q^t=q^{\^t}=q^4,\;
q^{y'}=q^4=q^{y''}
>
.
}
$$
Here,
\-
the first isomorphism is
\itemitem{--}
adding new generators
$y'$ and $y''$ and relators expressing the new generators via the old ones,
\itemitem{--}
and adding some consequences to the relator set;
relations (1) and (2) follows from the other relations:
$q^{y'}=q^{yt}=q^t=q^4$ (and similarly for (2));
\-
the second isomorphism removes some consequences from the
relator set;
relations (3) and (4) follows from the other relations:
$q^y=q^{y't^{-1}}=q^{4t^{-1}}=q$ (and similarly for (4));
\-
the third isomorphism removes the generators
$y$ and $\^y$ and relations (5) and (6)
that express these generators via the other ones.

\enditem
The kernel $\~H$ of the homomorphism from
$\BS(4,16)=\pres<a,b|a^{4b}=a^{16}>$ onto $\Z_4$,
mapping $a$ to 1 and $b$ to 0 has index four in $\BS(4,16)$.
If we take the set $T=\{1,a^{-1},a^{-2},a^{-3}\}$
as a Schreier transversal, then we obtain the following
generating set for $\~H$:
$$
\bigl\{\tau x(\={\tau x})^{-1}\;|\;x\in\{a,b\}\bigr\}
\setminus\1
=
\Bigl\{q=c^4,\;b,\;b'=b^a,\;b''=b^{a^2},\;b'''=b^{a^3}\Bigr\}
$$
(where $\=v\in T$ is the only element of $T\cap\~Hv$).
The defining relators for $\~H$ are words of the form
$\tau r\tau^{-1}$
(expressed via the generators of $\~H$),
where $r$ is a defining relator of the initial group
and $\tau\in T$.
Thus, we obtain
$$
\~H=\pres<q,b,b',b'',b'''|q^b=q^{b'}=q^{b''}=q^{b'''}=q^4>
$$
Comparing
this with the presentations of $H$,
we see that $H\iso\~H$, i.e., $\pi_1(L)$ is
commensurable with $\BS(4,16)$, which is not commensurable with
$\BS(2,4)=\pi_1(K)$ [CKZ21] (see also [Fo24]).  Thus, there is no finite
common covering for complexes $K$~and~$L$ (because the fundamental group
of a finite covering of a complex is isomorphic to a finite-index subgroup
of the fundamental group of this complex).

\s 4.
Why the universal coverings of $K$ and $L$ are isomorphic?

The universal covering of both complexes is
the \emph{Baumslag--Solitar complex} $X_{2,4}$ from [Fo24],
which has the following structure.

Consider the regular directed tree $T$ with all vertices of
out-degree two and in-degree four.
Let us assign
two numbers
$\gamma(e),\delta(e)\in\{0,1\}$
to each edge $e$ in such a way that
\-
each vertex serves as the head to
edges with all four values of the pair $\bigl(\gamma(e),\delta(e)\bigr)$,
\-
and each vertex serves as the tail to
edges with the both values of
$\gamma(e)$
and with the both values of $\delta(e)$.

\enditem
The Baumslag--Solitar complex $X_{2,4}$ is (homeomorphic to) the Cartesian
product $T\times\R$ with the following cell structure:
\-
the vertex set is $V(X_{2,4})=\{v_i\;|\;v\in V(T),\;i\in\Z\}$;
\-
the edges are the following:
\itemitem{--}
the edges $\epsilon_{v,i}$ (where $v\in V(T)$ and $i\in\Z$)
go from $v_i$ to $v_{i+1}$;
\itemitem{--}
the edges $e_i$ (where $e\in E(T)$ and $i\in\Z$),
go from $v_i$ to $w_{2i+\gamma(e)}$,
where $v$ and $w$ are the tail and head of $e$, respectively;
\-
the two-dimensional cells
$D_{e,i}$
are attached along the cycles of the form
$
e_i^{-1}
\epsilon_{v,i}
e_{i+1}
\epsilon_{w,2i+1+\gamma(e)}^{-1}\epsilon_{w,2i+\gamma(e)}^{-1}
$,
where $e\in E(T)$, $i\in\Z$, and $v$ and $w$ are the tail and head of~$e$,
respectively (Fig. 2).

\goodbreak
\bigskip
\centerline{\input 2.PIC}
\nobreak%
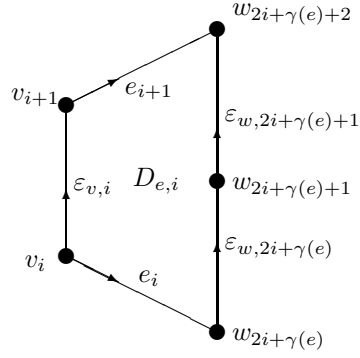
\centerline{Fig. 2}%
\goodbreak
\bigskip

\noindent
The covering $X_{2,4}\to K$ looks naturally:
\-
all vertices are mapped to the unique vertex of $K$;
\-
the edges are mapped by the following rule:
$\epsilon_{v,i}\mapsto c$;
$
e_i\mapsto\cases{
d,& if $i$ even,
\cr
z,& if $i$ odd;
}
$
\-
then the 2-cells $D_{e,i}$
are mapped to one of two 2-cells of $K$
depending on the parity of $i$.

\enditem
The covering $X_{2,4}\to L$ is the following:
\-
all vertices of $T$ are coloured black and white
in such a way that edges join vertices of different colours;
a vertex~$v_i$ is mapped to the vertex of $L$
of the same colour as $v$ (in what follows,
we consider a vertex $v_i$ coloured with the same colour as $v$);
\-
on the edges the mapping acts as follows:
\newline
$\epsilon_{u,i}\mapsto\cases{
c_\bullet& if $u$ is black,
\cr
c_\circ,& if $u$ is white,
}
$
and
$
e_i\mapsto\cases{
y,& if the tail of $e$ is black and $\delta(e)=0$,
\cr
z,& if the tail of $e$ is black and $\delta(e)=1$ is odd,
\cr
t,& if the tail of $e$ is white and $i+\delta(e)$ is even,
\cr
t_1,& if the tail of $e$ is white and $i+\delta(e)$ is odd,
}
$
\newline
(on Figure 3, $p,p',v,v',w,w',u,u'$ are vertices of~$T$,
and the pairs of numbers denote $\bigl(\gamma(e),\delta(e)\bigr)$
for the corresponding edge of~$T$);

\def\caption{\vbox{\small
\centerline{%
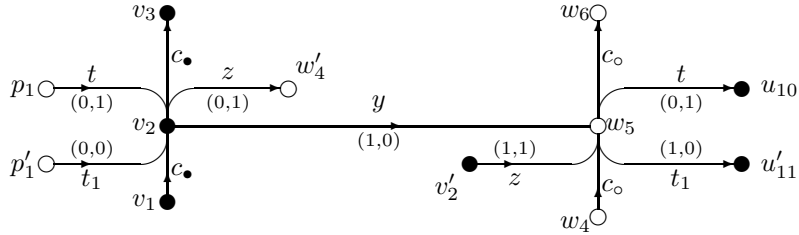
A part of the 1-skeleton of $X_{2,4}$
}
}}

\goodbreak
\bigskip
\centerline{\input 3.PIC}
\nobreak%
\centerline{Fig. 3}%
\goodbreak
\bigskip

\-
then the two-dimensional cell $D_{e,i}$
maps to
$
\cases{
A,& if the tail of $e$ is black, and $\delta(e)=0$,
\cr
B,& if the tail of $e$ is black, and $\delta(e)=1$,
\cr
D,& if the tail of $e$ is white, and $i+\delta(e)$ is even,
\cr
C,& if the tail of $e$ is white, and $i+\delta(e)$ is odd.
}
$

\enditem
It easy to verify that this map is a covering.
Figure 4 shows a neighbourhood of a (black) vertex $v$ in the tree $T$,
a neighbourhood of a corresponding vertex $v_2$ in the 1-skeleton
of $X_{2,4}$, and ten two-cells
($D_{e,0}$,
$D_{e,1}$,
$D_{\alpha,1}$,
$D_{\alpha,2}$,
$D_{f,0}$,
$D_{f,1}$,
$D_{\beta,1}$,
$D_{\beta,2}$,
$D_{g,0}$,
and
$D_{h,0}$),
adjacent to $v_2$; the framed letters denotes the images of edges and
two-cells under the covering $X_{2,4}\to L$.

\vfil\break

\goodbreak
\bigskip
\centerline{\input 4.PIC}
\nobreak%
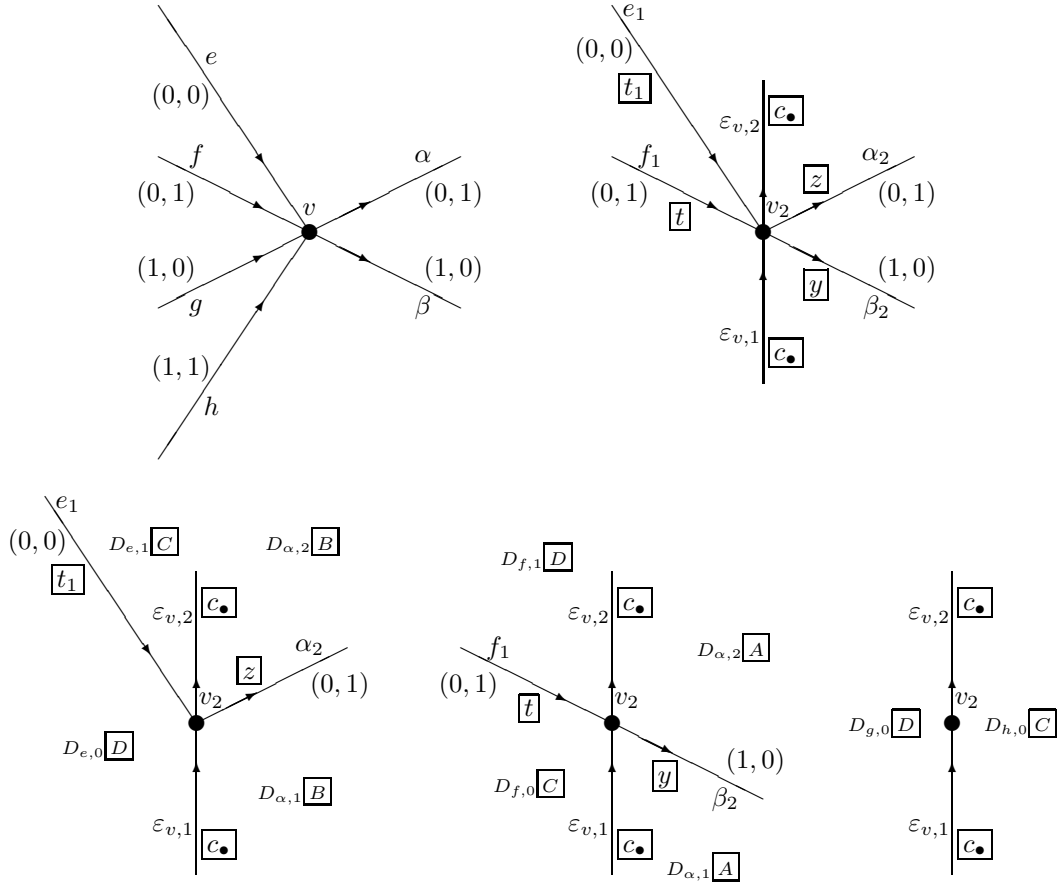
\centerline{Fig. 4}%
\goodbreak
\bigskip

\smallskip
\noindent
A brief description of this construction looks as follows:
\-
the complex $X_{2,4}$ is obtained from the regular tree $T$
by the following way:
\itemitem{--}
each vertex $v$ of $T$ corresponds to a line
$l_v$ in $X_{2,4}$;
this line consists of vertices $\{v_i\;|\;i\in\Z\}$ connected by edges
$\{\epsilon_{v,i}\;|\;i\in\Z\}$;
\itemitem{--}
each edge $e$ of $T$ corresponds to a band
$b_e$ in $X_{2,4}$ whose boundary
consists of two lines corresponding to ends of $e$;
this band $b_e$ is tessellated by edges $\{e_i\;|\;i\in\Z\}$
connecting the boundary lines;
\-
under the covering map $X_{2,4}\to L$, each band $b_e$ covers
\itemitem{--}
the annulus formed by the cell
$A$ in $L$
(this happens where $e$
is directed from the black vertex
to the white one
and
$\delta(e)=0$),
\itemitem{--}
or
the annulus formed by the cell
$B$ in $L$
(this happens where $e$
is directed from the black vertex
to the white one
and
$\delta(e)=1$),
\itemitem{--}
or
the annulus formed by the cells
$C$ and $D$ in $L$
(this happens where $e$ is directed from the white vertex to the black
one);
the mapping to $L$ on two such bands sharing a common boundary line are
shifted with respect to each other (because of different values of
$\delta(e)$).

\noindent
This complete the proof of the theorem.

\References

[BaK90]
H. Bass, R. Kulkarni,
Uniform tree lattices,
J. Amer. Math. Soc., 3:4 (1990), 843-902.

[BrS22]
M. Bridson, S. Shepherd,
Leighton's theorem: extensions, limitations, and quasitrees,
Algebraic and Geometric Topology, 22:2 (2022), 881-917.
\arXiv 2009.04305

[CKZ21]
M. Casals-Ruiz, I. Kazachkov, A. Zakharov,
Commensurability of Baumslag--Solitar groups,
Indiana Univ. Math. J., 70:6 (2021), 2527-2555.
\arXiv:1910.02117

[DK23]
N. S. Dergacheva, A. A. Klyachko,
Small non-Leighton two-complexes,
Math. Proc. Cambridge Philos. Soc., 174:2 (2023), 385-391.
\arXiv 2108.01398

[DK24]
N. S. Dergacheva, A. A. Klyachko,
Tiny non-Leighton complexes,
arXiv:2403.09803\thinspace.

[Fo24]
M. Forester,
Incommensurable lattices in Baumslag-Solitar complexes,
J. London Math. Soc., 109:3 (2024), e12879.
\arXiv:2207.14703

[JaW09]
D. Janzen, D. T. Wise,
A smallest irreducible lattice in the product of trees,
Algebraic and Geometric Topology, 9:4 (2009), 2191-2201.

[Lei82]
F. T. Leighton,
Finite common coverings of graphs,
J. Combin. Theory, Series B, 33:3 (1982), 231-238.


[LS77]
R. Lyndon, P. Schupp,
Combinatorial group theory.
Springer, 1977.

[MSW03]
L. Mosher, M. Sageev, K. Whyte,
Quasi-actions on trees I. Bounded valence,
Annals of Mathematics 158:1 (2003), 115-164.
\arXiv:math/0010136

[Neu10]
W. D. Neumann,
On Leighton's graph covering theorem,
Groups, Geometry, and Dynamics, 4:4 (2010), 863-872.
\arXiv 0906.2496

[SGW19]
S. Shepherd, G. Gardam,  D. J. Woodhouse,
Two generalisations of Leighton's Theorem,
arXiv:1908.00830.


[Sch27]
O. Schreier,
Die Untergruppen der freien Gruppen,
Abhandlungen aus dem Mathematischen Seminar der Universit\"at Hamburg,
5:1 (1927), 161-183.


[Wis96]
D. T. Wise,
Non-positively curved squared complexes:
Aperiodic tilings and non-residually finite groups.
PhD Thesis, Princeton University, 1996.

[Wis07]
D. T. Wise,
Complete square complexes,
Commentarii Mathematici Helvetici, 82:4 (2007), 683-724.

[Woo21]
D. Woodhouse,
Revisiting Leighton's theorem with the Haar measure,
Math. Proc. Cambridge Philos. Soc., 170:3 (2021), 615-623.
\arXiv 1806.08196


\end